\input amstex 
\documentstyle {amsppt} 
\vsize=8in \hsize=5.5in 
\topmatter  \title Special values of zeta-functions of regular schemes 
\endtitle
\author Stephen Lichtenbaum \endauthor

\endtopmatter

\document \magnification = \magstep 1 \baselineskip = 1.5 \baselineskip

\define \inl {\underleftarrow {lim}}
\define \ZZ {\Bbb Z}
\define \RR {\Bbb R}
\define \CC {\Bbb C}
\define \QQ {\Bbb Q}
\define \PP {\Bbb P}
\define \la {\lambda}

\S 0 Introduction.

In this paper we will give a conjectural formula (Conjecture 3.1)  for the special values $\zeta^*(X, r)$ of the scheme zeta-function of a regular scheme $X$ projective and flat of dimension d (so relative dimension $d-1$) over Spec $\ZZ$ at a rational integer $r$, in terms of singular, de Rham, and Weil-\`etale motivic cohomology, valid up to sign and powers of $2$. 

We can factor the map from $X$ to Spec $\ZZ$ uniquely through Spec $O_K$, where $K$ is a number field and the generic fiber of $X$ over Spec $O_K$ is a smooth connected algebraic variety over $K$.   We will construct complexes made up of variants of these cohomology groups, and the conjectured formula will give the special value as a product of Euler characteristics of these complexes, equipped with suitable integral structures.   We will prove that, if $d \leq 2$ this conjecture is compatible with Serre's conjectured functional equation for the zeta-function, and if $d > 2$ this compatibility is true modulo two previously existing conjectures which are true in dimension $\leq 2$.

We will discuss how this related to previous work on this  subject.  Beilinson, building on a previous conjecture of Deligne, gave a special values conjecture ([RSS]) for $r \leq 0$ up to a rational number.  More specifically, if one writes the scheme zeta function in the usual way as the product $\prod_{j=0}^{2d-2}L_j(X, s)^{(-1)^{j+1}}$, Beilinson gave a formula for the special value $L_j(X, r)$ up to a rational number.  Bloch and Kato ([BK]) gave a formula up to sign  for this special value when the weight $j-2r \leq -3$.  Fontaine and Perrin-Riou ([Fo]) gave such a conjectured formula  for all integers $r$.  Flach and Morin ([FM1]) gave a conjectured formula for special values of the zeta-function taking into account real primes. and more conceptual than that of Fontaine and Perrin-Riou.  In a very recent preprint [FM2], they also giver a proof of the compatibility of their conjecture with the functional equation, under essentially the same hypotheses as in the present paper. Fontaine asserts that the compatibility of his conjecture would follow from a local conjecture which still remains far from  being proven.  In this paper, as in ([FM1]), we only consider the zeta-function and not the associated L-functions.  We believe that the L-function conjectures are probably not completely correct, basically because of torsion phenomena in cohomology, which necessitate correction terms in analogous formulas for special values of zeta-functions of varieties over finite fields. We note that the wonderful formula of Bloch, Kato, and T. Saito, which plays an extremely important role in the proof of compatibility is only valid for Euler characteristics, and not for individual cohomology groups, which forces the restriction to zeta-functions. We also remark that throughout this paper we work with the scheme zeta-function of $X$. 

Recently Niranjan Ramachandran and I [(LR]) have shown that if $X$ is an arithmetic surface and $r=1$, rhe main conjecture of ths paper is equivalent to the conjecture of Birch and Swinnerton-Dyer for the Jacobian of $X$.  A similar result has been proved by Flach and Siebert ([FS]) for the main conjecture of [FM1].

The equivalence of the conjecture in this article with the  conjecture of Flach-Morin, and of either or both of these conjectures with the original conjecture of Fontaine and Perrin-Riou remains open, and is a very interesting and important question. 

There  are two basic approaches to zeta-function conjectures:  One (the Tamagawa approach) involves writing the formula as a product of local formulas (one for each prime p), and then using the product formula to show that although the individual factors may depend on choices, (possibly of a differential) the product does not.   This approach is used by Tate in his Bourbaki talk [T] on the conjecture of Birch and Swinnerton-Dyer for abelian varieties, by Fontaine and Perrin-Riou([Fo]), and by Flach-Morin. ([FM1]) . 

The other just involves just working with the infinite primes  and for example, choosing a particularly good differential.  This approach was used by the author in his conjectures on the Dedekind zeta function ([Li3]), by Silverman in stating the conjecture of Birch and Swinnerton-Dyer for elliptic curves in his book  ([Si]) on elliptic curves, and in this paper.   The basic idea of this paper is that by making the infinite prime part of the formula of Fontaine and Perrin -Riou more precise, we can dispense with the detailed local analysis.  

We should make it clear that, even to state our conjecture, we have to assume the validity of other previous conjectures.

First, we need the conjecture, which is very far from being proved, that the zeta-function of $X$, which converges for $Re (s)  > d$,can be meromorphically continued to the entire plane, so we can talk about $\zeta^*(X, r)$ for $ r < d$. 

Second, we assume that the \`etale motivic cohomology groups that we will define are finitely generated.

Third, we assume that the Beilinson regulator maps and various Arakelov intersection pairings induce isomorphisms on complex vector spaces. 

For the proof of compatibility, we need the theorem that  the groups $H^{2r+1}_{et}(X, \ZZ(r))$   and $H^{(2(d-r) +1}_{et}(X, \ZZ(d-r))$  
are finite and Pontriagin dual to each other.  This has been proved by Flach and Morin [FM1] if $d \leq 2$, and under some restrictions in the general case..

We also need, for the full Bloch-Kato-Saito theorem, resolution of singularities for arithmetic schemes.

\medskip

Let $X$ be a scheme of finite type over Spec $\ZZ$, and Krull dimension $d$.   Let $x$ denote a closed point of $X$, and let $N(x)$ be the order of the residue field $\kappa(x)$.   
Recall that $x$ is closed if and only if $\kappa(x)$ is finite.  Let $G$ be any meromorphic function  on $\CC$, let $r$ be a rational integer, and let $a_r$ be the order of the zero of $G(s)$ at $s = r$.  Let $G^*(r)$ be the limit as $s$ approaches $r$ of $G(s)(s-r)^{-a_r}$.  If $G$ is a zeta-function, $G^*(r)$ is referred to as a special value of $G$.

\proclaim {Definition 0.1} [Se1]: The zeta function $\zeta(X,s)$ of $X$ is defined to be $\prod_x (1-N(x)^{-s})^{-1}$, where $x$ runs over the closed points of $X$. \endproclaim

The product defining $\zeta (X, s)$ is well known to converge for $Re (s) > d$ and is conjectured to have a meromorphic continuation to the entire plane.  We will tacitly assume this conjecture in what follows.  It is further conjectured ([Se2]. [Bl3]) that there exists a $\Gamma$-factor $\Gamma (X, s)$ and a positive rational number $A$, the conductor, such that if we  let $\phi(X,s)= A^{s/2}\zeta(X,s) \Gamma(X,s)$ then $\phi(X, s)$ satisfies the functional equation 
$\phi(X, s) = \pm \phi(X, d-s)$.

Now let $X$ be regular, and projective and flat over Spec $\ZZ$.  The basic idea behind our conjectured formula is to start with Fontaine's ``Deligne-Beilinson'' conjectures [Fo], which give the special values up to a rational number in terms of determinants of maps of  complex vector spaces with given rational structures.  These complex vector spaces come from singular and de Rham cohomology, and from Weil-\'etale motivic cohomology. We replace the rational structures by integral structures, and take determinants with respect to these.  The singular cohomology of course has a natural integral structure, and the Weil-\'etale groups conjecturally do also.  We define an  integral structure on the de Rham groups by using derived exterior powers.  We should note that these derived exterior powers have an important role to play even in the number ring case ($d=1$).

  We also introduce the orders of naturally occurring finite cohomology groups into the picture.  Finally we replace the period maps in Fontaine's picture by ``modified'' period maps, where we divide by special values of the Gamma function.

  Our conjectural formula expresses the special values of the zeta-function in terms of the product of Euler characteristics of exact sequences of complex vector spaces with integral structures.  The complex vector spaces will be derived from singular cohomology, de Rham cohomology, and Weil-\`etale motivic cohomology.   For the exact formula, see Section 3.  The maps between them arise from Beilinson's conjectures, Arakelov height pairings, and periods.

The plan of this paper is as follows:

Section 1:  Integral structures and Euler characteristics.   

Section 2:  The groups and maps involved in the conjecture,

Section 3:  Statement of the conjecture.

Section 4:  The Euler characteristic of the period map.

Section 5: Serre's functional equation and $\Gamma$- function identities. 

Section 6:   Compatibility of the conjecture with the functional equation.

Section 7.  The case of number rings.

As we move along we will explain how our definitions of groups and maps relate to those of Fontaine and Perrin-Riou [FP].

We would like to thank Spencer Bloch, Matthias Flach, Thomas Geisser, Baptiste Morin, Niranjan Ramachandran, and Takeshi Saito for many helpful conversations.  

\S.1.  Integral structures and Euler characteristics.

Let $V_0, V_1 \dots V_n$ be finite-dimensional complex vector spaces, and 

$$ V_* = 0 \to V_0 \to V_1 \to \dots \to V_n \to 0$$

be an exact sequence.

Let $B_i$ be a lattice spanned by a basis for the vector space $V_i$.  Let $\Lambda V$ denote the highest exterior power of $V$ and $\Lambda B$ denote the highest exterior power of $B$ .  The alternating tensor product of the $\Lambda V_i's$ is canonically isomorphic to $\CC$ and the alternating  tensor product of the $\Lambda  B_i's$ is isomorphic to $\ZZ$.  The natural inclusions of $B_i$ in $V_i$ induce a map from $\ZZ$ to $\CC$ and the determinant $det (V_*, B_*)$ in $\CC^*/\pm 1$ of the pair $(V_*, B_*)$ is defined to be the image of a generator of $\ZZ$ in $\CC$.

\proclaim {Definition 1.1} Let $V_*$ be a sequence of finite-dimensional complex vector spaces. An integral structure on $V_*$ is a  sequence of pairs  $(A_*, a_*)$ where $A_*$ is a lattice in $V_*$ and $a_*$ is a positive rational number.   \endproclaim

An example of an integral structure on a finite-dimensional complex vector space $V$ comes from a finitely generated abelian group $M$ with a homomorphism from $M$ to $V$ whose image is a lattice $M_0$ in $V$ and whose kernel is the torsion subgroup $M_{tor}$  of $M$. The integral structure is then $(M_0, |M_{tor}|)$.

\proclaim{Definition 1.2} An integral structure $(V_*, A_*) $ is torsion-free if each $a_j$  is equal to $1$. \endproclaim

\proclaim {Definition 1.3} Let $ (A_*,a_*)$ be an integral structure on the finite exact complex $V_*$.  We define the Euler characteristic of $(A_*, a_*)$ to be det $(V_*,A_*) \prod a_j^{(-1)^{j+1}}$ \endproclaim

\proclaim {Definition 1.4} Let $(A_*, a_*)$ be an integral structure on $V_*$,  $( B_*, b_*)$ be an  integral structure on $W_*$, and $\phi_*$ be a map of complexes from $V_*$ to $W_*$ such that $\phi_j$ is an isomorphism for all $j$.  Let $det_j$ be the determinant of $\phi_j$ with respect to the lattices $ {A_j}$ and ${B_j}$, and let $\chi(\phi_j) = det_j b_j/a_j$.  Define the Euler characteristic $\chi(\phi)$to be $\prod_{j=0}^n (-1)^{j+1}\chi(\phi_j)$. \endproclaim

\proclaim {Definition 1.5} Let $\tilde{A} = ( A_*, a_*)$ be an integral structure on $V_*$. Let $V^\vee_*$ be the sequence of dual vector spaces to $V_*$. We define the dual integral structure $\tilde{A}^\vee$ on $ V^\vee _*$ to be  $(A^\vee_*, a_*^{-1})$, where $A^\vee_*$ is the lattice dual to $A_*$.: \endproclaim

\proclaim {Proposition 1.6} The Euler characteristic of the dual integral structure is equal to either to the Euler characteristic of the original integral structure or to its inverse, depending on whether $n$ is odd or even.\endproclaim

Proof.  Well-known

\S 2.  The groups and maps  involved in the conjecture.

\bigskip

\S 2.1 Weil-\'etale motivic cohomology..

\bigskip

As before,  Let $X$ be a regular scheme, projective and flat over Spec $\ZZ$  Let $X_0$ be the fiber of $X$ over Spec $\QQ$, and let $K$ be the algebraic closure of $\QQ$ in the function field of $X$, Let $O_K$ be the ring of integers in $K$. We may regard $X$ as a scheme projective and flat over Spec $O_K$.

We will first define Weil-\`etale motivic cohomology groups and then discuss their relation to the groups defined by Fontaine and Perrin-Riou ([Fo] and [Fl]). 

Let $r$ be an integer, and $j$ a non-negative integer.   We would like to define a Weil-\'etale site and complexes of sheaves  $\ZZ (r) $ on this site whose cohomology groups $H^j_W(X, \ZZ(r))$ would be Weil-\`etale motivic cohomology, but unfortunately we do not know how to do this.  Instead, for $j \leq 2r$ we define $H^j_W(X, \ZZ(r))$ to  be the hypercohomology groups  $H^j_{et}(X, \ZZ(r))$, where $\ZZ(r)$ denotes Bloch's higher Chow group complex sheafified for the \'etale topology ([Bl1].[Le]). Sometimes these groups are referred to as \'etale motivic cohomology.  
For $j \geq 2r +1$, we define $H^j_W( X, \ZZ (r))$ to be $h^j(RHom (R \Gamma_{et} (X, \ZZ(d-r)) , \ZZ[-2d-1])$, so we have the exact sequence 

$$ 0 \to Ext^1(H_{et}^{2d +2-j}  (X, \ZZ(d-r)), \ZZ) \to H^j_W(X, \ZZ(r))\to Hom (H_{et}^{2d +1 -j}(X. \ZZ(d-r)), \ZZ) \to 0 \leqno (2.1.1)$$ 

If we had our hypothetical Weil-\'etale site, with section functor denoted by $\Gamma_W$, this would follow, up to 2-torsion,  from a duality theorem which asserted that $R\Gamma_W(X, \ZZ(d-r))$ was isomorphic to $RHom (R\Gamma_W(X, \ZZ(r)), \ZZ[-2d-1])$.  The analogue of this theorem, assuming the usual conjectures, is true for Weil-\'etale cohomology in the geometric case, as shown in [Ge]).  We note here that in [FM1], Flach and Morin have constructed  such a complex of abelian groups, which satisfies this duality theorem assuming that standard finiteness conjectures hold.

The group $H^{2r}_W(X, \ZZ(r))$  is by definition $H^{2r}_{et}(X, \ZZ(r))$, and by standard arguments this agrees with the group $H^{2r}_{Zar}(X, \ZZ(r))$ of codimension $r$ cycles on $X$ modulo rational equivalence  after tensoring with $\QQ$.  Hence  there is a cycle map $\phi$ from $H^{2r}_W(X, \ZZ(r))$ to singular cohomology with rational coefficients.  Let $H^{2r}_W(X, \ZZ(r))^1$ denote Ker $\phi$ (cycles homologous to zero) and $H^{2r}_W(X, \ZZ(r))^2$ denote Image $\phi$ (cycles modulo homological equivalence.).

We have the exact sequence 

$$ 0 \to H^{2(d-r)}_W(X, \ZZ(d-r))^1 \to H^{2(d-r)}_W(X, \ZZ(d-r))  \to H^{2(d-r)}_W(X, \ZZ(d-r))^2 \to 0 $$

  \proclaim{Conjecture 2.1.1} The groups $H^j_{et}(X, \ZZ(r))$ are finitely generated for $j \leq 2r +1$, and finite for $j= 2r+1$. \endproclaim
                                   
This implies

\proclaim {Conjecture 2.1.2} The cohomology groups $H^j_W(X, \ZZ(r))$ are finitely generated for all j. \endproclaim

Assuming the truth of this conjecture, we give the complex vector space $H^j_W(X, \ZZ(r))_{\CC}$ the standard integral structure $H^j_W(X, \ZZ(r))$.

We also need 

\proclaim {Conjecture 2.1.3} The finite groups  $H^{2r+1}_{et}(X, \ZZ(r))$ and $H^{2(d-r)+1}_{et}(X, \ZZ(d-r))$ are Pontriagin duals. \endproclaim

 Flach and Morin show in  [FM1, Proposition 3.4] that this follws from Conjecture 2.1.2 for $d \leq 2$ and is  under some restrictions in the general case.

\bigskip

\S 2.2 Singular and de Rham cohomology.

We also will have need of singular cohomology groups.   Let $X_{\CC} = X \times_ \ZZ  \CC$.  
  Complex conjugation $c$ acts on $H^j_B(X_\CC, \ZZ)$ via the natural action of conjugation on $\CC$.  If $r$ is even (resp. odd), let $\tilde H_B^j(X, \CC(r))^+$ and  $\tilde H_B^j(X, \ZZ(r))^+$ be the set of elements $y$ in $H_B^j(X_\CC, \CC)$ and . $H_B^j(X_\CC, \ZZ)$  such that $c(y) = y$ (resp. $c(y) = -y$). We define $H_B^j(X, \CC(r)) $  (resp.$ H_B^j(X,\CC(r))^+$) to be $H_B^j(X_{\CC}, \CC)$ (resp. $\tilde H_B(X_\CC, \CC)^+)$ and its standard integral structure given by mapping $H_B^j(X_{\CC}, \ZZ)$  (resp. $\tilde H_B^j(X_\CC, \ZZ)^+$ ) to $H_B^j(X_{\CC}, \CC)$ (resp $\tilde H_B(X_\CC, \CC)^+) $ via the natural map followed by multiplication by $(2\pi i)^r$.

Let $\Omega = \Omega_{X_{\CC}/\CC}$. The de Rham cohomology group $H_{DR}^j(X_{\CC}, \CC)$ has the Hodge decomposition 

$$\prod _{i+k=j}H^i(X_{\CC}, \Lambda^k \Omega )$$ 

which gives rise to the Hodge filtration $G_m = \prod_{k\geq m} H^{j-k}(X, \Lambda^k \Omega)$. 

$H_{DR}^j(X_{\CC}, \CC(r))$ is defined to be $H_{DR}^j(X_{\CC}, \CC)$ but with the Hodge filtration $F$ given by $F_m = G_{m+r}$. If $M$ is $H^j(X, \ZZ(r))$ we define $t_M = t_{j,r}$ to be $ H^j_{DR}(X_{\CC} ,\CC(r))/F_0  =$

$=\prod_{k <r}H^{j-k} (X_\CC, \Lambda^k \Omega_{X_\CC/\CC}) = \prod_\sigma \prod_{k < r}H^{j-k}(X_\sigma, \Lambda^k \Omega_{X_\sigma/\CC})$.  Here $\sigma$ runs through all embeddings of the number field $K$ into $\CC$. and $X_\sigma = X \times_{O_K} \CC$ where the map from $O_K$ to $\CC$ is induced by $\sigma$. The standard integral structure on $ t_M$ is given by

 $\prod_\sigma \prod _{k < r} H^{j-k}(X, \lambda^k \Omega _{X/{O_K}})$, where $\lambda^k$ denotes the $k$th derived exterior power.  (See the appendix for a discussion of derived exterior powers).

 \bigskip

\S 2.3 Maps between cohomology groups.  

Let $M = M_{j,r}$ be the motive $H^j(X, \ZZ(r))$.  Let $M^* = H^{2d-2-j}(X, \ZZ(d-1-r))$, and $N = M^*(1) = H^{2d-2-j}(X, d-r))$. The classical period map $\alpha_{j,0}$ maps $H_B(M_{j,0})_\CC= H_B^j(X_\CC, \CC)$ to $H_{DR}(M)_\CC = H^j_{DR}(X_\CC, \ZZ_\CC$, which induces a map $\alpha_{j,r} = (2\pi i)^r\alpha_{j,0}$)  from $H_B(M)_{\CC}^+ = H^j_B(X, \ZZ(r))^+_{\CC}$ to $(t_M)_{\CC}$. We now define a new map $\gamma_M$  (which we call the enhanced period map) as follows:  $H_{DR}(M)$ has a decreasing Hodge filtration $F_{q'}(M)$.  Let $H^{q'} = F_{q'}/F_{q'+1}$.  Let $h_{q'}$ be the dimension of $H_{q'}$. Then $H_{DR}(M)$ has the direct sum Hodge decomposition $\coprod H^{q'}$.  We decompose $\alpha_M= \alpha_{j,r}$ into the direct sum of the maps $\alpha^{q'}(M)$ where $\alpha^{q'}$ is the map $\alpha_M$ followed by the projection onto $H^{q'}$.  Let $\Gamma$ be the usual gamma -function.  Recall that the weight $w(M)$ of $M$ is equal to $j-2r$,  Now let $\gamma^{q'}(M)$ be $ \Gamma^*(-w(M) +q')) \alpha^{q'}(M)$ and let $\gamma_{j,r}$ be the isomorphism $ \tilde \gamma_M = \coprod_{q'} \gamma^{q'}(M)$.  Since $h_{q'} = h(p,q)$, where $p+q=j$ and $q'=q-r$, the determinant of $\gamma_{j,r}$ is equal to the determinant of $\alpha_{j,r}$ multiplied by $\prod_{q'} \Gamma^*(-w(M) +q')^{h_{q'}}$ which is equal to $\prod _p\Gamma^*(r-p)^{h(p,q)}$, where $p+q=j$ and the product is over all $p$ between $0$ and $j$.

Consider the following diagram of exact sequences:

$$
\CD
0 @>>> (H_B(M)^+ )_{\CC}@>i>> (H_B(M))_{\CC}@>p>> (H_B(M)^-)_{\CC} @ >>> 0 \\
@.@VVV@V\tilde \gamma_M VV@VVV \\
0@>>> (F_0(M))_{\CC} @ >j>> (H_{DR}(M))_{\CC} @ >q>>( t_M)_{\CC} @ >>> 0 \\
\endCD
$$

Let $\gamma_M = q \circ \tilde \gamma_M \circ i $.  Let $\beta_M = p \circ \tilde \gamma_M^{-1} \circ j$. Diagram -chasing immediately shows that $\tilde \gamma_M$ induces isomorphisms from Ker $\gamma_M$ to 
Ker $\beta_M$ and  from Coker $\gamma_M$ to Coker $\beta_M$. 

\proclaim {Proposition 2.3.1}  The exact sequence of complex vector spaces 

$$ 0 \to Ker \gamma_M \to (H_B(M)^+ )_\CC \to (t_M)_\CC \to Coker \gamma_M \to 0 \leqno 2.3.1(M)$$ is dual to the exact sequence 

$$ 0 \to Ker \beta_N \to( F_0(N))_\CC  \to (H_B^-(N))_\CC \to Coker \beta_N \to 0 \leqno 2.3.2(N)$$. \endproclaim

Proof.  $F_0(N)_\CC$ is the Serre dual of $(t_M)_\CC$.  $H_B(N)^-$ may be identified with $H_B(M^*)^+$, which is the Poincar\`e dual of $H_B(M)^+$.  $H_B^{2d-2}(X, \CC(d-1)) $ may be canonically identified with $H_{DR}^{2d-2}(X, \CC)$  .  Poincar\`e duality is compatible with Serre duality, which implies the proposition.

We from now on choose an arbitrary basis for Ker $\gamma_M$, the basis for Ker $\beta_M$ induced by the isomorphism between Ker $\gamma_M$ and Ker $\beta_M$, the basis for Coker $\gamma_N$ induced by the above duality, and the  basis for Coker $\beta_N$ induced by the isomorphism between Coker $\gamma_N$ and Coker $\beta_N$   We will use these integral structures on the various kernels and cokernels,

If $A$ is a finitely generated abelian group, Let $A_{tor}$ denote the torsion subgroup of $A$ and $\bar {A}$ denote $A/A_{tor}$.  

If $\phi: A \to B$ is a homomorphism of finitely generated abelian groups, let $\bar {\phi}$ be the induced homomorphism from $\bar {A}$ to $\bar {B}$ and let $\phi_{tor}$ be the induced homomorphism from $A_{tor}$ to $B_{tor}$.

\proclaim {Lemma 2.3.2} .  Let $ 0 \to A \buildrel f \over  \to B \buildrel g \over \to C \to 0$ be an exact sequence of finitely generated abelian groups.  There is a natural isomorphism from Ker $\bar {g}/ Im \bar {f}$ to Coker $g_{tor}$, and the determinant of $0 \to \bar A \to \bar B \to \bar C \to 0$ is equal to the Euler characteristic $|B_{tor}|/|A_{tor}||C_{tor}|$ of $0 \to A_{tor} \to B_{tor} \to C_{tor} \to 0$.  \endproclaim

Proof. Exercise

\proclaim {Proposition 2.3.3} $ \chi(2.3.2(M) ) = \chi (\tilde{\gamma_M}) \chi( 2.3.1(M))$ \endproclaim

Proof.  Let $\Lambda$ denote highest exterior power. Let $A_1^* = H_B(M)^+$, let $A_2^* = H_B(M)$ and let $A_3^*= H_B(M)^-$.   Let $A_j$ be a generator of $\Lambda (\bar{A_j^*})$.  Let $B_1^* = F_0(M)$, $B_2^* = H_{DR}(M)$, and $B_3^* = t_M$.  Let $B_j$ be a generator of $\Lambda (\bar{B_j^*})$.  Let $C_1$ be a generator of any lattice $L_1$ in $\Lambda (Ker (\beta_M) )= \Lambda(Ker (\gamma_M))$ and $C_2$ be a generator of any lattice $L_2$ in $\Lambda (Coker (\beta_M)) = \Lambda(Coker (\gamma_M))$.

Let $A = A_2/A_1A_3$ and $B = B_2/B_1B_3$.  Let $a_j = |(A_j^*)_{tor}|$ and $b_j = |(B_j^*)_{tor}|$.  Let $\chi(A_{tor}) = a_1a_3/a_2$ and $\chi(B_{tor}) = b_1b_3/b_2 $.  $\chi_{tor} (2.3.1(M)) = b_3/a_1$, and $\chi_{tor} (2.3.2(M) = a_3/b_1$.

We have $\chi_{tor}(2.3.1(M))/\chi_{tor}(2.3.2(M)) = b_1b_3/a_1a_3$, and $\chi(2.3.1(M))/\chi(2.3.2(M))=$

$$ = (det (2.3.1(M))/det(2.3.2(M))) (\chi_{tor}(2.3.2(M))/\chi_{tor}(2.3.1(M)))= $$  
 $$= B_1B_3a_1a_3/A_1A_3b_1b_3= = (B_2 det(B)/A_2 det(A)) (\chi(A_{tor}a_2/\chi(B_{tor}b_2$$
 
 Lemma 2.3.2 implies that this equals $B_2a_2/b_2A_2$, which is $\chi(\bar{\gamma}_M)^{-1}$.

Beilinson defines Chern class maps from the algebraic K-theory groups to Deligne cohomology. Let $\gamma_{j,r} =\gamma_M$ where $M=h^j(X,\ZZ(r))$.  In our language, Beilinson's map becomes a map $c_{j,r}$ (for $j\leq 2r-2$) from $H^{j+1}_W(X, \ZZ(r))_{\CC}$ to  Coker $\gamma_{j,r}$.  Let $N$ be the motive $M^*(1) = h^{2d-2-j}(X,\ZZ (d-r))$.  Since Ker $\gamma_N$ may be identified with the dual of Coker $\gamma_{2d-2-j,d-r}$,  we also have the dual $b_{j,r}$ (for $j\geq 2d-2r)$ of Beilinson's Chern class map which maps Ker $\gamma_{2d-2-j,d-r}$ to $H^{2d-j}_W(X, \ZZ(d-r))_{\CC}$.  We have:

\proclaim{Conjecture 2.3.4} (Beilinson)  If  $j \leq 2r-3$ the maps $c_{j,r}$ and $b_{j,r}$ are isomorphisms.  \endproclaim

Recall that by definition $H^{2r+1}_W(X, \ZZ(r))_{\CC}$ is dual to $H^{2(d-r)}_{et}(X, \ZZ(d-r))_{\CC}$ which is the same as $H^{2(d-r)}_{Zar}(X, \ZZ(d-r))_{\CC}$, i.e. codimension (d-r) cycles on $X$ modulo rational equivalence.  Recall that $H^{2(d-r)}_{et}(X, \ZZ(d-r))^1$ denote cycles homologically equivalent to zero, and $H^{2(d-r)}_{et}(X, \ZZ(d-r))^2$ denote cycles modulo homological equivalence.  Let $H^{2r+1}_W(X, \ZZ(r))^1_{\CC}$ be the dual of $H^{2d-r}_{et}(X, \ZZ(d-r))^1_{\CC}$

\proclaim{Conjecture 2.3.5} (Beilinson)  There is an exact sequence

$$  0 \to H^{2r-1}_W(X, \ZZ(r))_{\CC} \buildrel c_{2r-2,r} \over \to Coker  (\gamma_M) \to H^{2r+1}_W(X, \ZZ(r))_{\CC}^2 \to 0 $$ 

with $M = H^{2r-2}(X, \ZZ(r))$.

\endproclaim

This is a slightly different but more natural variant of Beilinson's original conjecture, and it is implicitly used by Fontaine ([Fo]).

\proclaim {Conjecture 2.3.6} There is an exact sequence

$$ 0 \to H^{2r}_W(X, \ZZ(r))^2_{\CC} \to Ker (\gamma_N) \to H^{2r+2}_W(X, \ZZ(r))_{\CC} \to 0 $$

\endproclaim

This is the dual of conjecture 2.3.5, with $M$ replaced by $N=M^*(1)= H^{2d-2r}(X, \ZZ(d-r))$ 

\proclaim {Conjecture 2.3.7} There is an isomorphism $(H^{2r+1}_W(X, \ZZ(r))^1)_{\CC}$ to $(H^{2r}_W(X, \ZZ(r))^1)_{\CC}$

\endproclaim

This is the non-degeneracy of the Arakelov intersection pairing restricted to finite cycles homologous to zero, where it is independent of metrics.  

\bigskip \bigskip

\S.3 The statement of the conjecture.

\bigskip

We would like to first explain the relationship between Weil-\`etale motivic cohomology groups and the groups which occur in Fontaine's Deligne-Beilinson conjecture.  We look at the motive $M =H^j(X, \ZZ (r))$,  Recall that $N = M^*(1)$ is $H^{2d-2-j}(X, \ZZ(d-r))$.  

Fontaine starts with a projective non-singular algebraic variety $X_0$ over Spec $\QQ$.  He chooses a regular model $X$ for $X_0$ projective and flat over Spec $\ZZ$.  He conjectures that the following six-term sequence is always exact:  

$$ 0 \to H^0_f(M)_{\CC}  \to Ker (\gamma_M)  \to H^1_c(M)_{\CC} \to H^1_f(M)_{\CC} \to Coker (\gamma_M) \to H^2_c(M)_{\CC} \to 0 $$.

 If $j \leq 2r-2$, Fontaine's $H^1_f(M)$ is our $H^{j+1}_W(X, \ZZ(r))_{\QQ} $, (We are using motivic cohomology instead of algebraic K-theory, but these two groups agree after tensoring with $\QQ$). (Actually, Fontaine's group is the image of $K(X)$ in $K(X_0)$, but we conjecture that the natural map is always injective.)  

If $j = 2r-1$, Fontaine's  $H^1_f(M)$ is the group of codimension $r$ cycles on $X_0$ homologically equivalent to zero, tensored with $\QQ$,

 If $j \geq 2r$, $H^1_f(M) = 0$

Fontaine's $H^0_f(M)$ is zero unless $j = 2r$, in which case it equals the group of codimension $r$ cycles on $X_0$ modulo homological equivalence, tensored with $\QQ$.

Fontaine's $H^i_c(M)$ is the $\QQ$-dual of $H^{2-i}_f(N)$ for $i=1,2$.

For each $ j $ and $r$ with  $j \leq 2r-3$ or $j\geq 2r+1$, we will define a sequence of integral structures $A(j,r)$ and $A'(j,r)$.
$A(j,r)$ is given by:

$$ (j \leq 2r-3)  \quad \quad      c_{j,r}:  H^{j+1}_W(X, \ZZ(r))_{\CC} \to Coker (\gamma_{j,r}) $$

while $A'(j,r)$  is given by:
 $$ (j \geq 2r +1) \quad \quad b_{j,r}:  Ker (\gamma_{j,r}) \to H^{j+2}_W(X, \ZZ(r))_{\CC}  $$
 
If $1<r<d$ we define an integral structure $C(r)$ given by: $C(r)=$
 
 $$ 0 \to H^{2r-1}_W(X, \ZZ(r))_{\CC} \buildrel c_{2r-2,r} \over \to Coker (\gamma_{2r-2,r}) \to H^{2r+1}_W(X, \ZZ(r)) _{\CC} \buildrel e_r \over \to $$
$$ \to H^{2r}_W(X, \ZZ(r))_{\CC} \to Ker (\gamma_{2r,r}) \buildrel b_{2r,r} \over \to H^{2r+2}_W(X. \ZZ(r))_{\CC}  \to 0 $$

Here $e_r$ is induced by the Arakelov intersection pairing.
We give these vector spaces the standard integral structures previously defined in \S. 2.2.

 Conjectures 2.3.4, 2.3.5, 2.3.6 and  2.3.7 imply that these sequences are exact.
 
 We give degrees to the terms of these complexes by requiring that $Ker (\gamma_M)$ has even degree and $Coker (\gamma_M) $ has odd degree.  (These sequences are all truncations of modified versions of Fontaine's six-term sequence in [Fo] and this convention makes the degrees agree)
  
  Finally we define exact sequences $B(j, r)_\CC$  given  for all $j$ and $r$ by
 
 $$ 0 \to ( Ker(\gamma_M) \to (H^j_B(X, \ZZ(r))^+)_{\CC} \buildrel \gamma_M  \over  \to t_{j,r} \to Coker(\gamma_M) \to 0. $$
 
 We put $Ker(\gamma_M)$ in degree zero.
 
 The integral structures on the cohomology groups here are induced by the standard integral structures defined in \S 2.2.

 Let $\chi_{A,C} (X, r) = \chi(C(r))  \prod _{j=0}^{min(2d-1,2r-3)} (\chi(A(j,r)) ^{(-1)^{j}}  \prod_{j=max (0,2r+1)}^{2d-1} \chi(A'(j,r))^{(-1)^{j}}$
 
 Let $\chi _B(X,r) = \prod_{j=0}^{\infty}\chi(B(j,r))^{(-1)^{j}} $  Let $\chi(X,r) = \chi_{A,C}(X,r) \chi_B(X,r)^{(-1)}$
  
  \proclaim {Conjecture 3.1} Give all groups in the  above exact sequences their standard integral structures. 
    $$\zeta^*(X,r)=  \chi (X,r) $$ up to sign and powers of 2.(Note that $A(j,r)$ and $B(j,r)$ are torsion for $j\geq 2d-1$, and zero for $j$ large).\endproclaim
  
  If $j \neq 2r-1$ each of the terms Ker $(\gamma_M)$ and Coker $(\gamma_M)$ occurs exactly twice in the conjecture with degrees of opposite parity, so the conjecture is independent of the choice of integral structure.  If $j = 2r-1$, Ker $(\gamma_M)$ and Coker $(\gamma_M)$ are both zero.  
  
  \proclaim {Proposition 3.2}  a)If $0 \leq j\leq 2r-3$, then   $2(d-r)+1 \leq 2d-2-j \leq 2d-2$  and the integral structures $A'(2d-2-j, d-r)_{tf}$ and $A(j,r)_{tf}$ are dual..   Hence $det(A(j,r))= det(A'(2d-2-j, d-r))$  
 
  b) The integral structures $C(r)_{tf}$ and $C(d-r)_{tf}$ are dual.  Hence $det (C(r)) = det (C(d-r))$.
   \endproclaim
   
   \proclaim {Proposition 3.3} $\chi_{A,C}(X, r) = \chi_{A,C}(X, d-r)$ \endproclaim.
   
   Proof.  $\chi(C(r)) = det (\bar{C(r)})/ tor (C(r))$. $\bar {C(r)}$ is dual to $\bar{C(d-r)}$, so $det (C(r) = det C(d-r)$.  $\chi(A(j,r) = det (A(j,r))/ tor (A(j,r)$.  $\chi(A'(j,r)) = det (A'(j,r))/ tor (A'(j,r)).$ $\bar{A(j,r)}$ is dual to $\bar {A'(2d-2=j, d-r)}$ so $det(A(j,r)) = det (A'(2d-2-j, d-r))$ and $det(A'(j,r) = det A(2d-2-j, d-r))$.  
   
   On the other hand
   
  $$tor C(r)\prod_{j=0}^{2r-3}tor(A(j,r))^{(-1)^j} \prod_{j=2r+1}^{2d-1} tor (A'(j,r)^{(-1)^j} = \prod_{j=1}^{2d+1} |H^j(X,\ZZ(r))_{tor}|^{(-1)^j} $$
  which is equal by duality to 
  
  $$ \prod_{j=1}^{2d+1}|H^j(X, \ZZ(d-r))_{tor}|^{(-1)^j}$$
  
  which is equal to 
  
  $$tor(C(d-r))\prod _{j=0}^{2d-2r-3} tor (A(j, d-r))^{(-1)^j} \prod _{j=2d-2r+1}^{2d-1} tor (A'(j, d-r))^{(-1)^j}$$
   \bigskip
   \bigskip

\S 4  The Euler characteristic of the period map.

\bigskip

Let $K$ be the integral closure of $\QQ$ in the function field of $X$.  Then we can view $X$ as a scheme over $S = $ Spec $O_K$, and $X_{\CC}$ is canonically isomorphic to $\prod_\sigma X \times _S$ Spec $\CC$, where the product is taken over all embeddings $\sigma$ of $K$ in $\CC$.  Let $\kappa(v)$ be the residue field of the closed point $v$, and let $X_v= X\times _S $ Spec $\kappa(v)$

Recall that $A$ is the positive rational number which appears in the conjectured functional equation $A^{s/2} \Gamma(X, s) \zeta(X,s) = \pm A^{(d-s)/2}\Gamma(X, d-s) \zeta(X, d-s)$ for $\zeta(X,s)$.  Let $\omega_{X/S}$ be the relative canonical class of $X$ over $S$.

\proclaim {Definition 4.1} $A'_v = (\Delta_X. \Delta_X)_v= (-1)^d {c_d}^X_{X_v}(\Omega_{X_v/O_v}) \in CH_0(X_{\kappa(v)})$.  $A' =  (\Delta_X. \Delta_X)_S =\prod_v A'_v$.\endproclaim

(This definition is taken from [KS].  We will not use it directly in what follows.  What we need is stated in Theorems 4.2 and 4,3)

Here the product is taken over all closed points $v$ of $S$ such that $X$ is not smooth over $S$ at $v$. Note that the Chern class ${c_d}^X_{X_v}(\Omega)$ is equal to 1 if X is smooth over $S$ at $v$.

\proclaim {Theorem 4.2} (Bloch-Kato-Saito) If $d \leq 3, A' = A$.  If strong resolution of singularities holds for schemes of finite type over Spec $\ZZ$, then $A'=A$ in general. \endproclaim

Proof.  This is the main result (Theorem 6.2.3) of [KS](the conductor formula of Bloch).

  Let $O_v = O_K$ localized at $v$. Let $X_v = X \times_{O_K} O_v$.

\proclaim {Theorem 4.3} The Euler characteristic of the cone $C_{m, v}$ of the map from 

$R\Gamma(X_v, \lambda^m \Omega_{X_v/O_v})$ to $R\Gamma(X_v,RHom (\lambda^{d-1-m} \Omega_{X_v/O_v} , \omega_{X_v/O_v}))$ induced by Serre duality is equal to $(A')_v^{(-1)^{m+1}}$.  \endproclaim

Proof.  This is Corollary 4.9 of [Sa].

\proclaim{Theorem 4.4} The Euler characteristic of the classical period map $\alpha$ from de Rham to singular cohomology of $X_\CC$ with respect to the canonical integral structures is 

$(A')^{d/2} (2\pi i)^{\chi(X_{\CC})(d-1)/2}$.  \endproclaim

This theorem will follow from the following propositions.

\proclaim {Proposition 4.5} Let $H^j$ be $H^j_{DR}(X_{\CC}, \CC) = \prod_\sigma H^j_{DR} (X_\sigma. \CC)$ with its canonical integral structure.  Let $
G^j$ be $H^j_{DR}(X_{\CC}, \CC)$  with the integral structure given by 

$\prod_\sigma \prod _{k=0}^j Ext^{2d-2-j+k}_{X_\sigma}(\lambda^{d-1-k}\Omega_{X/{O_K}}, \omega_{X/O_K})$.  Then the Euler characteristic $\chi_{I}$ of the identity map from $H^*$ to $G^*$ with respect to the given integral structures is $(A')^d$. \endproclaim

Proof. Let $\Omega = \Omega_{X_\sigma/\CC}$ and $\Omega_K = \Omega_{X/{O_K}}$,
By Serre duality, $H^{j-k}(X_\sigma, \Lambda^k \Omega) $ is canonically isomorphic  to $Ext^{k-j+1-d}_{X_\sigma} (\Lambda^{d-1-k}\Omega, \omega_{\CC})$.  So $ A_{j,k} = H^{j-k}(X, \lambda^k\Omega_K)$ and .$B_{j,k} = Ext_X^{j-k}(\lambda^{d-1-k}\Omega_K, \omega)$ give two different integral structures on $H^{j-k}(X_\sigma, \Lambda^k\Omega)$.  

Let $d_{j,k}$ be the determinant of the identity map with respect to these two integral structures.  (This is independent of the choice of basis for the integral structures, up to sign.) Then Theorem 4.3 asserts that $\prod_j (d_{j,k})^{(-1)^{j-k}}. \chi((A_{j,k})_{tor})^{-1} \chi ((B_{j,k})_{tor}) =  (A')_\sigma $ if $k$ is even and is equal to $(A')_\sigma^{-1}$ if $k$ is odd.  By Serre duality, $\chi((A_{j,k})_{tor})= (\chi((B_{j,k})_{tor})^{-1}$.

We conclude easily from this that $\prod_k \prod_j (d_{j,k})^{(-1)^j} (\chi (\lambda^k(\Omega))_{tor}))^{-2} = (A')^d$ and the product on the left is the Euler characteristic of the identity map. 

Let $Q_{DR}$ denote the de Rham cup product, and $Q_B$ denote the cup product in singular cohomology.   Let $\alpha^j: H^j_{DR}(X_{\CC}, \CC) \rightarrow  H^j_B(X_{\CC}, \CC)$ be the classical period map,  We know that $Q_{DR}$ and $Q_B$ are compatible, i.e. $Q_B(\alpha^j(a) , \alpha^k(b)) = \alpha^{j+k}(Q_{DR}(a, b))$

\proclaim {Proposition 4.6} Let $  [v_{i,j} ]$ be a basis of $H^j_{DR}(X_{\CC}, \CC)$.  Let $[u_{i,j}]$ be a basis of $H^j_B(X_\CC, \CC)$ coming from a basis of $H^j(X_{\CC}, \ZZ )$ modulo torsion. .  Let $E_j$ be the matrix of $\alpha^j$ with respect to the bases $[ v_{i,j}]$ and $[ u_{i,j}]$. Then

$$ (\prod_j (det E^j)^{(-1)^j})^2 = \prod_j det (Q_{DR}( v_{i,j}. v_{k, 2d-2-j}) )^{(-1)^j} (2\pi i)^{\chi(X_{\CC})(d-1)}$$
\endproclaim

Proof.    Fix $j$.  Since on the top-dimensional cohomology group  $H_{DR}^{2d-2}(X_{\CC}, \CC)$, $\phi(1_{DR})= (2\pi i)^{d-1} 1_B$ , we have  $\phi(Q_{DR}(v_{i,j} , v_{k,2d-2-j} )) = Q_B(\phi(v_{i,j}), \phi(v_{k, 2d-2-j}))$ which implies 

$(2\pi i)^{(d-1) B_j} det (Q_{DR}(v_{i,j}, v_{k, 2d-2-j})) = det Q_B(\phi^j (v_{i,j}) , \phi^{2d-2-j} (v_{k, 2d-2-j})) = $

$= det (E^j) det (E^{2d-2-j}) det(Q_B(u_{i,j}, u_{k, 2d-2-j}))= \pm det(E^j) det(E^{2d-2-j})$.  Multiplying by $(-1)^j$ and taking products over $j$, we obtain the propositiion.                                                                                                                                                                                                                                                                                                                                                                                                                                                                                                                                                                                                                                                                                                                                                                                                                                                                                                                                                                                                                                                                                                                                                                                                                                                                                             

\proclaim {Proposition 4.7} Let $d_j = \prod_k d_{j,k}$ be the determinant of the identity map with respect to the two given integral structures on $H^j_{DR}(X_\CC,\CC)$. Then a) $d_j = det ( Q_{DR}(x_{i,j}, x_{k, 2d-2-j}))$

b) $ \prod d_j^{(-1)^j} = \prod det ( (Q_{DR}( x_{i,j},  x_{k,2d-2-j})^{(-1)^j} )$ \endproclaim

Proof.  By the compatibility of cup product and Serre duality, $det ((Q_{DR}x_{i,j}, y_{k,2d-2-j})$  is equal to  $1$.  Then 4.7a follows from the definition of $d_j$, and of course implies 4.7b..

Proof of Theorem 4.4:

Recall that by definition and Serre duality $\chi(\alpha) = \prod_j (det E^j)^{(-1)^j} \chi((H_{DR})_{tor})^{(-2)}$ 

and $\chi_{I} = \prod d_j^{(-1)^j} \chi((H_{DR})_{tor})^{-1}$. Substituting $x$ for $v$ in Proposition 4.6  and using Proposition 4.5 we obtain

 $$(A')_\sigma^d (\chi(H^*_{DR}(X)_{tor})^2 = ((\prod_j det (Q_{DR}([x_{i,j}], [x_{k, 2d-2-j}])^{(-1)^j)})$$ 

Proposition 4.6 implies

$$(A')_\sigma^d (\chi(H^*_{DR}(X)_{tor})^2 = (2\pi i)^{-\chi(X_{\CC})(d-1)}(\prod det(\alpha^j)^{(-1)^j})^{-2} $$

which implies

$$(2\pi i)^{\chi(X_{\CC}) (d-1)/2} (A')_\sigma^{d/2} = \det (\alpha^*) \chi ((H^*_{DR}(X))_{tor})^{(-1)} $$

By Poincar\`e duality, $\chi(H^*_B(X_{\CC}, \ZZ)_{tor} ) = 1$, so

$$\det (\alpha^*) \chi((H^*_{DR}(X))_{tor})^{(-1)} \chi((H^*_B(X_{\CC}, \ZZ)_{tor}) = (2\pi i)^{\chi (X_{\CC})(d-1)/2 }(A')_\sigma^{d/2} $$

$$\chi(\alpha) = (2 \pi i)^{(\chi(X_{\CC})(d-1)/2 } (A')_\sigma^{d/2}$$

Taking the product over $\sigma$, we get Theorem 4.4.

\proclaim {Corollary 4.8}  The Euler characteristic $\chi_r$ of the classical period map $\alpha_r$ from 

 $H^*_{DR}(X_\CC, \CC (r))$ to $H^*_B(X, \CC, \CC(r))$ is $(A')^{d/2} (2\pi i)^{\chi(X_\CC)((d-1)/2 -r)}$.  \endproclaim
 
 This follows from the definition of twisting by $r$.

Recall that  $\Gamma_{r,j} = \prod_{p+q=j}\Gamma^*(r-p)^{h(p,q)}$.  Let $\Gamma_r = \prod \Gamma_{r,j}^{(-1)^{j+1}}$

\proclaim {Corollary 4.9} The Euler characteristic  $\chi(\gamma_r)$ of the enhanced period map $\gamma_r$ is 

$$\Gamma_r( A')^{d/2} (2\pi i)^{\chi(X_\CC)(((d-1)/2 -r)}$$\endproclaim

Proof:  By the remarks at the beginning of section 2.3, $\chi(\gamma_r)$ is equal to $\chi(\alpha_r)$ multiplied by $\prod_{j=0}^{2d-2} \Gamma_{r,j}^{(-1)^{j+!}} = \Gamma_r$

\bigskip \bigskip

\S 5.1  Serre's functional equation and $\Gamma$-function identities

\bigskip

Let $X_0$ be a smooth projective algebraic variety of dimension $d -1$ over the number field $K$.  Let $j$ be a non-negative integer, and let $L^j(X_0,s)$ be the L-function attached by Serre in [Se2] to the $j$-th cohomology group of $X_0$,

Let $\sigma$ be an embedding of $K$ into $\CC$, and 
let $v$ be the place of $K$ induced by $\sigma$.  Let $K_v$ be the completion of $K$ at $v$.   Let $X_v = X_0 \times_K \CC$, where $\sigma$ maps $F$ into $\CC$, and let $\Omega =\Omega_{X_v/\CC}$. Recall that Hodge theory gives us a decomposition $H_{DR}^j(X_v)= \coprod H_v^{p.q}$, where the sum is taken over pairs $(p.q)$  such that $p+q =j$ and $ H_v^{p.q} = H^q(X_v, \Omega^p)$. Let $c$ be the automorphism of $X_v$ induced by complex conjugation acting on $\CC/K_v$.  Then if $j$ is even and equal to $2n$, $c $ acts as an involution on $H_v^{n,n}.$ Let $h_v(p,q)$ be the dimension of $H_v^{p,q}$.

Then $H_v^{n,n} = H_v^{(n,+)} \oplus H_v^{(n,-)}$, where

$$ H_v^{(n,+)} = \{x \in H_v^{n,n} , c (x) = (-1)^n x\} $$

$$H_v^{(n,-)} = \{ x \in H_v^{n,n}, c (x) = (-1)^{(n+1)} x\} $$

Let $h_v(n, +) = dim H_v^{(n,+)}$, and $h_v(n -) = dim H_v^{(n,-)}$.

Let $B^j_v$ be the rank of $H^j(X_v, \ZZ)$ and let $(B_v^j)^+$ be the rank of the subgroup of $H^j(X_v, \ZZ)$ left fixed by $c$.  Let $(B_v^j)^- = B^j_v -(B^j_v)^+$.   Note that if $j = 2n$, $(B_v^j)^+$ is equal to $\Sigma h_v(p,q) + h_v(n,+)$ if $n$ is even and is equal to $\Sigma h_v(p,q) + h_v(n,-)$ if $n$ is odd, where the sum is taken over all pairs $(p,q)$ where $p<q$ and $p+q =j$.  Let $(B_v^{j,r})^-$ be $(B_v^j)^-$ if $r$ is even, and $(B^j_v)^+$ if $r$ is odd.

Let $\Gamma_{\CC}(s) = (2\pi)^{-s}\Gamma(s)$.  Let $\Gamma_{\RR}(s) = \pi^{-s/2}\Gamma(s/2)$.  

Serre gives the functional equation $\phi^j(s) =\pm \phi(j+1-s)$, where $\phi(s) = L^j(s) A_j^{s/2} \Gamma^j(s)$, $A_j$ is a certain positive integer, and $\Gamma^j(s)$ is described as follows:

$\Gamma^j(s) = \prod _v \Gamma_v^j(s)$
where $\Gamma^j_v(s) = \prod_{p+q=j}(\Gamma_{\CC}(s-inf(p,q))^{h_v(p,q)}$ if $v$ is a complex place of $F$,

$\Gamma_v(s)= \Gamma_{\RR}(s-n)^{h_v(n,+)} \Gamma_{\RR}(s-n+1)^{h_v(n,-)}\prod _{p <j-p}\Gamma_{\CC}(s-p)^{h_v(p,j-p)}$, if $v$ is a real place of $F$.

We observe that it is an easy computation that $\Gamma_v^j(s) = \Gamma_v^{2d-2-j}(s+d-j-1)$, so that with our earlier observation that at least in the smooth case $L^j(s) = L^{2d-2-j}(s+d-j-1)$ we obtain that Serre's functional equation is equivalent to the functional equation $\phi^j(s) = \pm \phi^{2d-2-j}(d-s)$.

\proclaim {Theorem 5.1.1}Let $v$ be a real place of $K$. 

 If $j$ is even, $ (\Gamma^j_v)^*(r)/ (\Gamma^{2d-2-j}_v)^*(d-r) $ is equal up to sign and powers of 2 to
 $\prod_p (\Gamma^*(r-p))^{h_v(p,q)} \pi^{-{B_v^j(r - j/2) +(B_v^{j,r})^-}}$.  (The products run over $0 \leq p \leq B_v^j$).

If $j$ is odd, $ (\Gamma^j_v)^*(r)/ (\Gamma^{2d-2-j}_v)^*(d-r) $ is equal up to sign and powers of 2 to $\prod_p \Gamma^*(r-p)^{h_v(p,q)} \pi^{-B_v^j(r -(j+1)/2)}$ \endproclaim 

.  Proof.   We consider the case when $j$ is even. (The case when $j$ is odd is similar, but simpler). Let $j=2n$. Fix $v$ and let $q = j-p$. First look at terms where $p\neq q$.     Let $p' = d-1-p$ and $q' = d-1 - q$, so $p' + q' = 2d-2-j$.  We have 

$$\prod_{p<q}\Gamma_{\CC}^*(r-p)^{h_v(p,q)}/\prod_{p'<q'} \Gamma_{\CC}^*(d -r -p')^{h_v(p',q')}  =$$,

$$\prod_{p<q}\Gamma_{\CC}^*(r-p)^{h_v(p,q)}/\prod _{p>q} \Gamma_{\CC}^*(1-r+p)^{h_v(p,q)} $$

since $h_v(p,q) = h_v(p', q')$ by Serre duality.

By definition of $\Gamma_{\CC}$, this is equal to 

$$ \prod_{p<q} \Gamma^*(r-p)^{h_v(p,q)}/\prod_{p>q} \Gamma^*(1-r+p)^{h_v(p,q)} $$

multiplied by

 $$(2\pi) ^{-(\Sigma_{p<q}(h_v(p,q) (r-p) -\Sigma_{p>q}(h_v(p,q)(1-r+p) ) }$$

This product is then equal to 

$$ \pm \prod_{p\neq n} \Gamma^*(r-p)^{h_v(p,q)} (2\pi)^{-(\Sigma_{p,q} (h_v(p,q)( (r-p) - (1-r+j-p))}  $$

because of the relation $\Gamma^*(r) = \pm \Gamma^*(1-r)^{-1} $ for integral $r$ which follows immediately from the functional equation for the Gamma function.  We then obtain:

$$\pm  \prod_{p\neq n} \Gamma^*(r-p)^{h_v(p,q)} (2\pi)^{-(B_v^j-h_v(n,n)) (r - (j+1)/2)} \leqno (5.1.1)$$

We now look at the terms involving $n$ with $v$ still fixed.

We first observe that the functional equation for the Gamma function implies that

$\Gamma^*(a/2) \Gamma^*((2-a)/2)$ equals  $\pm \pi$ if $a$ is an odd integer and equals $\pm 1$ if $a$ is an even integer.

We compute:

$$\Gamma_{\RR}^*((r-n))^{h_v(n,+)} \Gamma_{\RR}^*((r-n+1))^{h_v(n,-)}$$
multiplied by
 $$\Gamma_{\RR}^*((d-r-(d-1-n))^{-h_v(n,+)} \Gamma_{\RR}^*(d+1-r-(d-1-n))^{-h_v(n,-)}$$

which we rewrite as 

$$(\Gamma^*((r-n)/2) (\Gamma^*(1-r+n)/2)^{-1})^{h_v(n,+)} \leqno (5.1.2)$$

multiplied by

$$ \Gamma^*((r-n+1)/2) (\Gamma^*((n+2-r)/2)^{-1} )^{h_v(n,-)} \leqno (5.1.3)$$

multiplied by 
$$\pi^{-( h_v(n,+)(2r-2n-1)/2 +h_v(n,-)( 2r-2n-1)/2 )} \leqno (5.1.4)$$

First assume that $r-n$ is odd.  By the functional equation,
 $$\Gamma^*((1-r+n)/2)^{-1} = \pm \Gamma^*(r-n+1)/2 \leqno (5.1.5)$$ 

 $$\Gamma^*((n+2-r)/2)^{-1} = \pm \pi^{-1} \Gamma^*((r-n)/2) \leqno (5.1.6)$$.

Now recall the duplication formula for the Gamma function;

$$\Gamma(2s) = 2^{2s-1}\Gamma(s) \Gamma(s + 1/2) /\sqrt \pi \leqno (5.1.7)$$

Then (5.1.2) becomes (using (5.1.5) and (5.1.7)) 

$$(\pm 2^{n-r+1} \Gamma^*(r-n) \sqrt \pi)^{h_v(n,+)} $$

and (5.1.3) becomes

 $$(\pm 2^{n-r} \Gamma^*(r-n) \sqrt \pi)^{-h_v(n,-)} $$

while (5.1.4) is

$$\pi^{-h_v(n,n)(r-(j+1)/2} \leqno (5.1.8)$$

So, up to sign and powers of 2, our product has become

$$ \Gamma^*(r-n)^{h_v(n,n)} \pi^{(h_v(n,n)(r- (j+1)/2) +h_v(n,+)/2 -h_v(n,-)/2} \leqno (5.1/9) $$

Multiplying  (5.1.1) by (5.1.9) we get 

$$\prod_p \Gamma^*(r-p)^{h_v(p,q)} \pi^{-(B^j_v(r -(j+1/2))) +h_v(n,+)/2 -h_v(n,-)/2) } \leqno (5.1.10) $$
which equals

$$ \prod_p\Gamma^*(r-p)^{h_v(p,q)} \pi^{-(B^j_v(r-j/2) +(B^j_v)^+} \leqno (5.1.11) $$

if $n$ is even ( so $r$ odd) and equals

$$\prod_p\Gamma^*(r-p)^{h_v(p,q)} \pi^{-B^j_v(r-j/2) +(B^j_v)^-} \leqno (5.1.12)$$

if $n$ is odd (so $r$ even).
since $(B_j -h_v(n,n)/2) +h_v(n,+)$ is equal to $(B^j_v)^+$ if $n$ is even and $(B^j_v)^-$ if $n$ is odd.

The proof for $r-n$ even  is identical, except for switching $h_v(n,+)$ and $h_v(n, -)$.

\proclaim {Theorem 5.1.2} Let $v$ be a complex place of $K$.  Then $(\Gamma^j_v)^* (r) /(\Gamma^{2d-2-j}_v)^*(d-r)$ is equal up to sign and powers of 2 to $(\prod_p\Gamma^*(r-p))^{2h(p,q)} \pi^{-B^j_v(2r - (j+1))} $ \endproclaim

Proof. Let $q =j-p$, $p'= d-1-p$ and $q'= d-1-q$. Note that $0 \leq, p, q, p', q', \leq d-1$

We write  $(\Gamma^j)^*_v(r)/(\Gamma^{2d-2-j})^*_v(d-r) $ = $$\prod_p \Gamma_{\CC}^*(r-Inf(p,q))^{h_v(p,q)}/\prod_{p'} \Gamma_{\CC}(d-r-Inf(p',q'))^{h_v(p',q')}$$ 

which is equal to

$$\prod_p\Gamma^*(r-Inf(p,q))^{(h_v(p,q)} /\prod_{p'}\Gamma^*(d-r -Inf (p',q'))^{h_v(p',q')} \leqno (5.1.13)$$

multiplied by $$(2\pi)^{-\Sigma_p h_v(p.q) (r- Inf(p,q)) -\Sigma_{p'}h_v(p',q')(d-r-Inf(p',q'))} \leqno (5.1.14) $$

Since $h_v(p',q') =h_v(p,q)$, the functional equation for the Gamma function transforms (5.1.13) into (up to sign) 

$$ \prod_p \Gamma^*(r-Inf(p,q))^{h_v(p,q)}   \prod_p \Gamma^*(r-Sup(p,q))^{h_v(p,q)}  \leqno (5.1.15) $$,

which in turn equals 

$$(\prod_p\Gamma^*(r-p)^{h_v(p,q)})^2 $$

Now (5.1.14) easily transforms into $ (\pi^{2r-j-1})^{-\Sigma_p h_v(p,q)}$, which becomes $\pi^{-(B^j_v(2r - (j+1))}$.

\bigskip

\S.6  Compatibility of the conjecture with the functional equation.

Starting with Serre's conjectured functional equation for cohomological L-functions described in the previous section, Bloch, Kato, and T. Saito  conclude that the following functional equation holds for the zeta-function of $X$;

\proclaim {Conjecture 6.1} Let $\phi (X,s) = \zeta (X, s) A^{s/2} \Gamma (X, s)$.  Then $\phi (X, s) = \pm \phi (X, d-s)$ \endproclaim

Here the constant $A = A(X)$ is obtained by taking the alternating product of the constants $A_j$ which occur in the conjectured functional equation for Serre's L-function $L_j$ and modifying it by terms coming from degenerate fibers.  It is still a positive rational number. 

The generic fiber $X_0$ of $X$ is a projective algebraic variety smooth over a number field $K = H^0(X_0, O_{X_0})$.  Let $\Gamma (X, s) = \prod_j \prod _{\sigma} \Gamma(X_{v(\sigma)}^j, s)^{(-1)^j}$, where $\Gamma_v^j$ was defined in the previous section. 
If we rewrite our conjecture as $\zeta^*(X, r) = \chi(X,r)$ , then what we want to show is that $\zeta^*(X, r)/\zeta^*(X, d-r) = \chi(X, r)/\chi(X, d-r)$, up to sign and powers of 2, where the left-hand side is computed by the functional equation.
\proclaim {Proposition 6.2}Conjecture 6.1 implies 

$$\zeta^*(X, r)/\zeta^*(X, d-r) = A^{d/2 -r}\prod_j (\prod_{p+q=j} {{\Gamma^*(r-p)}^{h(p,q)})}^{(-1)^j} \pi^{-\chi(X_{\CC}) (r-(d-1)/2) +\chi^-(X_{\CC}, \ZZ(r))}$$. \endproclaim

This is an immediate consequence of Theorems 5.1.1 and 5.1.2, remembering that $B_j = B_{2d-2-j}$.

We now wish to compute  $\chi(X, r)/\chi(X, d-r)$ and show that it agrees with the expression in Proposition 6.2.  We first recall that $\chi(X, r) = \chi_{A,C}(X,r) \chi_B(X, r)$ and that $\chi_{A,C}(X,r)/\chi_{A,C}(X,d-r) = 1$, by  Proposition 3.3.  

We now have to look at $\chi_B(X, r)/\chi_B(X, d-r)$.

\proclaim {Lemma 6.3} Let $\chi_{i,k}$ be the Euler characteristic of the identity map from the complex vector space $H^i(X_{\CC}, \Lambda^k \Omega_{X_\CC})$ with the  integral structure $H^i(X, \lambda^k \Omega_X)$ to the same vector space with the integral structure $\underline {RHom} (H^{d-1-i} (X,  \lambda^{d-1-k}\Omega_X),  \omega)$.  Then $\prod_i \chi_{i,k}^{(-1)^i} =( A')^{(-1)^k}$. \endproclaim

Proof. This follows immediately from Theorem 4.3.
\proclaim {Lemma 6.4} Let $\theta_{j,r} $ be the Euler characteristic of the identity map from the complex vector space $t_M = \coprod_{0\leq k < r} H^{j-k}(X_{\CC}, \Lambda^k \Omega_{X_{\CC}})$ with the integral structure $\coprod _{0\leq k  < r} H^{j-k}(X, \lambda^k \Omega_X )$  to the same vector space with the integral structure given by $$ \underline {RHom} (H^{d-1-j+k}(X, \lambda^{d-1-k} \Omega_X), \omega)$$.  Then $\prod_j \theta _{j,r}^{(-1)^j} = (A')^r $.\endproclaim 

Proof.  This follows immediately from Lemma 6.3

\proclaim {Lemma 6.5} Let $\eta_j$ be the Euler characteristic of the identity map from  $H^j_B(X_\CC ,\CC)^+$ (if $r$ is odd) or $H^j_B(X_\CC, \CC)^-$ (if $r$ is even) with the integral structure $H^j_B (X, \ZZ(r-1))^-$ 
 to the same vector space with the integral structure $H^j_B(X, \ZZ(r))^+ $  or $H^j_B(X, \ZZ(r))^-$.  Then  $\prod \eta_j ^{(-1)^j} = (2\pi i)^ {\chi^-(X_\CC, \ZZ(r))}$ \endproclaim
Proof.  This follows immediately from the definitions.

Let $\gamma(j,r)$ denote the integral structure $ 0 \to Ker \gamma_* \to H^j_B(X, \CC(r))^+ \to t(j,r)_\CC \to Coker \gamma_* \to 0 $

Let $\beta(j,r)$ denote the integral structure $ 0 \to Ker \beta_* \to {F_0(j,r)}_\CC   \to H_B^j(X, \CC(r))^-   \to Coker \beta_* \to 0 $

\proclaim {Proposition 6 .6} a) The integral structure $\beta (j,r)$ is dual to the integral structure $\delta(2d-2-j,d-r)$ given as follows:

$$ 0 \to Ker \gamma_* \to (H_B^{2d-2-j}(X, \ZZ(d-1-r))_\CC)^-  \to t(2d-2-j, d-r)_\CC \to Coker \gamma_* \to 0 $$.

where the singular cohomology group has its standard structure and t(2d-2-j, d-r) has the integral structure dual to the standard one on $F_0(j,r)$.

b) det $(\beta(j,r)) = det (\delta(2d-2-j, d-r))$ \endproclaim

Proof. Part a) follows immediately from the compatibility of Serre and Poincar\`e duality, and then b) follows immediately.

Let $\chi(\gamma(r)) = \prod_{j=0}^{2d-2} \chi(\gamma(j,r))^{(-1)^j}$, $\chi (\beta(r)) = \prod _{j=0}^{2d-2} \chi(\beta(j,r))^{(-1)^j}$. and $\chi(\epsilon (d-r)) = \prod_{j=0}^{2d-2} \chi( \epsilon (j, d-r))^{(-1)^j}$.

\proclaim {Corollary 6.7} $\chi(\beta(r)) = \chi(\delta(d-r))$ \endproclaim

Proof. This follows from the preceding proposition and a consideration of torsion, exactly as in Proposition 3.3.

Note that $\gamma(2d-2-j, d-r)$  is the integral structure $\delta (2d-2-j, d-r)$ except that the singular cohomology group now has the integral structure given by $H_B^{2d-2-j}(X, \ZZ(d-r)_\CC)^+$ and the tangent space has the standard integral structure rather than the one coming from duality.

\proclaim {Proposition 6.8} $\chi(\gamma(d-r)) = (A')^{d-r}   (2\pi i)^{  -\chi^-(X_\CC, \ZZ(r)) }  \chi (\delta(d-r))$ \endproclaim
Proof.  This follows immediately from Lemmas 6.4, 6.5 and the previous  remark.  Note also that $\chi^-(X_\CC, \ZZ(d-1-r)) = \chi^-(X_\CC, \ZZ(r))$ because Poincar\`e duality is compatible with complex conjugation.

\proclaim {Proposition 6.9}$ \chi(\gamma(r))/\chi(\beta(r)) = \Gamma_r(A')^{d/2} (2\pi i) ^{\chi(X_\CC) ((d-1)/2 -r) }$\endproclaim

Proof.  This follows from Proposition 2.3.3 and  Corollary 4.9

\proclaim {Proposition 6.10} $  \chi(\gamma(d-r) )/\chi(\gamma(r) ) = \Gamma_r(A')^{d/2-r}(2\pi i)^{\chi(X_\CC)((d-1)/2 -r) + \chi^-(X_\CC, \ZZ(r))}$ \endproclaim

Proof;  This follows immediately from Propositions 6.8, 6.9, and Corollary 6.7

\proclaim {Theorem 6.11}  $\chi_B(X, d-r) /\chi_B(X, r) = \Gamma_r(A')^{d/2-r}(2\pi i)^{\chi(X_\CC)(r-(d-1)/2 ) + \chi^-(X_\CC, \ZZ(r))}$ \endproclaim

 Proof.  $\chi_B(X, d-r) /\chi_B(X, r) =\chi( \gamma (d-r)) /\chi(\gamma (r)) $

Proposition 6.2 and Theorem 6.11 immediately imply the compatibility of our conjecture with the functional equation, if we replace $A$ by $A'$.

\proclaim {Theorem 6.12} If $d \leq 2$, our conjecture is compatible with the functional equation.  \endproclaim This follows from Theorem 6.11 and the main theorem of [KS].

\S 7.  The case of number rings

Let F be a number field with ring of integers $O_F$, class number $h$, number of roots of unity $w$, and discriminant $d_F$. Let $X = $ Spec $O_F$.  We will explain how our  conjecture for $X$ and $r$ reduces to standard theorems if $r=0$ or $r=1$ and well-known conjectures if $r <0$ or $r >1$.  

We begin with $r=0$.

We know that $H^j_{et}(X, \ZZ(1)) = 0$ if $j < 1$, $H^1_{et}(X, \ZZ(1)) = O_F^*$, $H^2_{et}(X, \ZZ(1))= Pic(X)$, $H^3(X, \ZZ(1)) = 0$ (up to 2-torsion) and $H^0_{et}(X, \ZZ(0)) = \ZZ$. 

 It immediately follows from the definitions that $H^0_W(X, \ZZ(0)) = \ZZ$, $H^1_W(X, \ZZ(0))= 0$,   and $H^3(X, \ZZ(0)) = \mu_F^\vee$, the dual of the roots of unity in $F$.  
 
 We also have the exact sequence $0 \to Pic(X)^\vee \to H^2_W(X, \ZZ(0)) \to Hom(O_F^*, \ZZ) \to 0 $.
 
 We also have $t_{j,0} = 0$ for all $j$.  It follows that $A(j,0)$ is always equal to $0$.
 
 We also see that $A'(j,0) = 0$ unless $j = 1$, when $A(1,0)$ reduces to $(\mu_F)^\vee$ in degree 2, so $\chi(A'(1,0)) = w$.
 
 $C(0)$ becomes $0 \to \CC \to \CC^{r_1+r_2} \to Hom (O_F^*, \CC) \to 0 $,  with the integral structure on the last term being $H^2_W(X, \ZZ(2)$ and the  second map being the dual of the classical regulator. So $\chi(C(0))= hR$, and $\chi(X, 0) = hR/w$.  As  is well-known, $\zeta^*(X,0) = -hR/w$.
 
 We now consider the case when $r =1$.
 
 $A(j,1)$ is easily seen to be zero for $j < -1$.
 
 The complex part of $\CC(1)$ is given by $0 \to O_F^*\otimes \CC  \to \CC^{r_1+r_2} \to \CC \to 0 $, with the last two terms getting the standard bases and the first term a basis coming from $O_F^*$, so  $det (\CC(1)$ is the classical regulator $R$. The Euler characteristic of $\CC(1)_{tor}$ is  $$|H^!_W(X, \ZZ(1))_{tor}|/|H^2(X,\ZZ(1))_{tor}| =w/h$$.  It follows that the Euler characteristic of $\CC(1)$ is $hR/w$.
 
 Since $t_{j,1} = 0$ for $j \geq 0$ , and $H^j_W(X, \ZZ(1)) = 0$ for $j \geq 5$, $A'(j,1) = 0$ for $j \geq 3$.
 
 Finally, $B(j.1)=0$ if $j \neq 0$, and $B(0,1)$ is given by 
 
 $$ 0 \to \CC^{r_2} \to O_F \otimes \CC \to \CC^{r_1 +r_2} \to 0 \leqno (7.1)$$. 
 
Note that we have the usual map $\theta$ mapping $O_F \otimes \CC$ to $\CC^{r_1+2r_2}$ by  sending $x \otimes 1$ to the collection of $\sigma(x) $ as $\sigma$ runs through the embeddings of $F$ in $\CC$.  The map from $\CC^{r_2}$ to $O_F \otimes \CC$ is given by the natural inclusion of $\CC^{r_2}$ in $\CC^{r_1 +2r_2}$ multiplied by $2\pi i$, followed by the inverse of $\theta$.  Since the determinant of $\theta$ with  respect to the usual bases is $\sqrt d_F$, we see that the determinant of (7,1) is equal to $(2\pi i)^{r_2}/\sqrt d_F$.  Hence the Euler characteristic $\chi(X, 1)$ is equal to $hR (2\pi i)^{r_2}/w\sqrt d_F$ which  is equal to the usual formula $hR(2\pi)^{r_2}2^{r_1}/w\sqrt |d_F|$ for $\zeta^*(X, 1)$ up to a power of 2.

Now let $r <0$.

The only non-zero groups $H^k_W(X, \ZZ(r))$ occur when $k =2$ or $k =3$, so $A(j,r)$ is equal to zero for all $j$ in the appropriate range, $C(r) =0$, and $A'(j,r) = 0$ unless $j=0$ or $j =1$. 

$B(j,r) = 0$ unless $j =0$.  $B(0,r)$ reduces to the isomorphism $Ker( b_{0,r}) \to H^2_W(X, \ZZ(r))^+$ so we give $Ker (b_{0,r})$ the integral structure induced from $H^2_W(X, \ZZ(r))^+$, and $\chi(B(0,r)) =1$.

$A'(0,r)_\CC$ is given by $Ker (b_{0,r}) \to H^2_W(X,\ZZ(r))_\CC$.  $H^2_W(X, \ZZ(r))_\CC$ is dual to $H^1_W(X, \ZZ(1-r))_\CC$ which is canonically isomorphic to $K_{1-2r}(O_F)_\CC$, and $Ker (b_{0,r})$ is canonically dual to $Coker (c_{0,1-r})$, where $c_{0,1-r}$ is the Beilinson regulator map. We conclude that $det(A'(0,r))$ is dual to the determinant $D_{1-r}$ of the Beilinson regulator map with respect to the natural bases coming from singular cohomology and K-theory.

Then $\chi(X,r)$ is equal to $D_{1-r} |H^2_W(X, \ZZ(r))_{tor}|/|H^3_W(X, \ZZ(r)|$, which is equal to $D_{1-r} |H^2_W(X, \ZZ(1-r)|/|H^1_W(X, \ZZ(1-r))_{tor}|$ which in turn is equal to $D_{1-r} |K_{-2r}(O_F)|/|K_{1-2r}(O_F)_{tor}|$, up to 2-torsion.

This is essentially what was conjectured in [Li3] to be $\zeta^*(X, r)$

Finally, let $r > 1$.  By definition, since $j$ has to be between $2r+1$ and $2d-1$, there is no contribution from $A'(j,r)$

$A(j,r)_\CC$ is the complex $H^{j+1}_W(X, \ZZ(r)) _\CC\to Coker (\gamma_{j,r})$, which is only non-zero when $j =0$, in which case the map is the Beilinson regulator $c_{0,r}$,  Since $j$ has to be either $0$ or $1$, the torsion Euler characteristic is $|H^2(X, \ZZ(r))|/|H^1(X, \ZZ(r))_{tor}$ which up to $2$-torsion is $K_{2r-2}(O_F)|/|K_{2r-1}(O_F)_{tor}|$.  So letting $R_r$ be the determinant of the Beilinson regulator map, the Euler characteristic $\chi(A,r))= (|K_{2r-2}(O_F)/|K_{2r-1}(O_F)_{tor}|)R_r$, up to 2-torsion.

Now $B(j,r)_\CC$ obviously is zero if $j \neq 0$.  Let $s(r) = r_2$ If $r$ is even and $s(r )= r_1+r_2$ if $r$ is odd.  Then by the same arguments as in the case when $r =1$ we have $det(B(0,r)_\CC = (2\pi i)^{rs(r)}\sqrt(d_F)$.

The integral structure on $t_{j,r} $ is given by $\coprod_{0\leq k < r}H^{j-k}(X, \lambda^k (\Omega))$  In the appendix, we compute that $\lambda^k(\Omega)$ is equal to $\Omega[k-1]$.  So $t_{j,r}= \coprod_{k<r} H^{j-k}(X,\Omega[k-1])$.  Since $\Omega$ only has cohomology in dimension $0$, we see that $t_{j,r} =0$ unless $j =1$, and the order of $H^0(X, \Omega)$ is $d_F$ so $\chi_{tor}(t_{j,r}) = d_F^{-r}$, and $\chi(B(j,r)) = (2\pi i)^{rs(r)} d_F^{(r-1/2)}$.

Putting everything together, we get that our conjecture says that $\zeta(X,r) = \chi(X,r)$, where  
$$\chi(X,r) = ((K_{2r-2}(O_F)|/|K_{2r-1}(O_F)_{tor}|)((2\pi i)^{rs(r)} d_F^{r-1/2}R_r $$

 up to 2-torsion, which is compatible with our conjecture for $s=1-r$ via the functional equation, which of course here is a well-known theorem.

Appendix:  Derived Exterior Powers

Let $\Cal A$ be an abelian category .  Let $S\Cal A$ denote the category of simplicial objects of $\Cal A$ and $C \Cal A$ denote the category of homological chain complexes of objects of $\Cal A$ ending in degree zero.  There are well-known functors $N: S\Cal A \to C\Cal A$ and $K: C\Cal A \to S\Cal A$ such that $NK$ is the identity and $KN$ is naturally equivalent to the identity functor.  $N$ and $K$ also preserve homotopies.  Let $\Lambda^k$ denote  $k$-th exterior power.  Let $X$ be a scheme and $\Cal A$ be the category of coherent locally free sheaves on $X$. If $Q_\cdot$ is in $S\Cal A$ with $Q_n$ a locally free sheaf on $X$  for all $n$. we define $\Lambda^k Q_{\cdot}$ to be $\Lambda^k (Q_n)$ in $S\Cal A$.

\proclaim {Proposition A.1}

Let $X$ be a regular scheme projective over Spec $\ZZ$.  Write $X$ as a closed subscheme of a projective space $P = (P^n)_{\ZZ}$ such that $I$ is the sheaf of ideals defining $X$ .  
Then the complex of locally free sheaves $C_{X,P} =I/I^2 \to \Omega_{P/\ZZ}$ defines an element in the derived category of locally free sheaves on $X$ which is independent of the choice of embedding of $X$ into $P$. \endproclaim

Proof.  If we have two different embeddings of $X$ in $P_1$ and $P_2$, take the Segr\'e embedding of $P_1 \times P_2$ in $P_3$, and compare successively the complexes defined by the embeddings into $P_1$ and $P_2$ with the product embedding into $P_3$.  (For details, see [LS]).

\proclaim {Definition A.2}  $\lambda^k \Omega_{X/\ZZ} = N\Lambda^k KC_{X,P}$. \endproclaim 

We see easily that this is independent of the choice of embedding.

We begin by recalling the following fact ([H], Exercise 5.16 (d)):

\proclaim  {Lemma A.3} Let $(X, O_X)$ be a ringed space, and
let $0 \to F' \to F  \to F'' \to 0 $ be an exact sequence of locally free sheaves of  $O_X$-modules.  Then there exists a finite filtration of $\Lambda^rF$:

$$ \Lambda^r F = G^0 \supseteq G^1   \supseteq \dots \supseteq G^r    \supseteq G^{r+1} = 0 $$

with quotients $G^p/G^{p+1} = \Lambda^pF' \otimes \Lambda^{r-p}F''$. \endproclaim

\proclaim{Theorem A.4} Let $(X,O_X)$ be a scheme such that every coherent sheaf on $X$ has a finite resolution by locally free sheaves.   Let $0\to F' \to F \to F'' \to 0$ be an exact sequence of coherent sheaves on $X$.  Let $r$ be a positive integer,  Then there exist objects in the derived category $G^0, G^1, \dots G^{r+1}$ and maps from $G^p$ to $G^{p+1}$ such that we have $G^0 = \lambda^r(F), G^{r+1} = 0$, and exact triangles $G^p \to G^{p+1} \to \lambda^p(F') \otimes ^L \lambda^{r-p}(F'') \to G^p[1]$. \endproclaim

Proof,  This is an easy corollary of the definition of derived exterior power and the previous lemma.

\proclaim {Theorem A.5} Let $A$ be a ring and $M$ an $A$-module.
a) $\lambda^0 M$ is canonically homotopic to the complex consisting only of $A$ in degree $0$.

b) $\lambda^1 M$ is canonically homotopic to the complex consisting only of $M$ in degree 0.

c) If $M$ has projective dimension $r$ then $\lambda^k M$ is represented by a complex of length $kr$.

\endproclaim

Proof: a) and b) are  obvious and c) follows immediately from a theorem of Dold and Puppe ([DP].

\proclaim {Theorem A.6} Let $X = $Spec $O_F$., and $\Omega = \Omega_{O_F/\ZZ}$.  Then $\lambda^k(\Omega)$ is isomorphic to $\Omega[k-1]$ \endproclaim 

Proof.  Let $D$ be the inverse different of $O_F$ over $\ZZ$.  We use induction on $k$, applying  Theorem A.4 to the exact sequence $0 \to O_F \to D \to \Omega \to 0$, and using that $\lambda^k(O_F) = \lambda^k(D) = 0$ for $k >1$.

References

[Bl1] Bloch , S. Algebraic cycles and higher K-theory, Adv, Math 61 (1986). 267-304

[Bl2] Bloch, S. Algebraic Cycles and the Beilinson Conjectures, Contemp. Math. Vol. 58 Part I, (1986), 65-78

[Bl3] Bloch, S. De Rham cohomology and conductors of curves.  Duke Math. J. 54 (1987) no.2 295-308

[BK] Bloch, S.; Kato, K.  L-functions and Tamagawa numbers of motives. The Grothendieck Festschrift, Vol. I, 333-400, Progr. Math., 86, Birkhauser Boston, Boston, MA, 1990.

[DP] Dold, A. and Puppe, D. Homologie nicht-aadditiven Funktoren.  Anwendungen.  Annales de l'Institut Fourier 11 (1961) 201-312

 [Fl] Flach, M. The equivariant Tamagawa number conjecture: a survey. With an appendix by C. Greither. Contemp. Math., 358, Stark's conjectures: recent work and new directions, 79-125, Amer. Math. Soc., Providence, RI, 2004.
 
 [FM] Flach, M. and Morin, B., Weil-\'etale cohomology and Zeta values of proper regular arithmetic schemes. Doc, Math 23 (2018) 1425-1560
 
 [FM2] Flach, M. and Morin, B. Compatibility of Special value conjectures with the functional equation of Zeta functions.  arXiv 2005.04829

[Fo] Fontaine, J-M.  Valeurs sp\'eciales des Fonctions-L des motifs, Sem. Bourbaki 34 (1991-1992), Exp. No. 751

[FP] Fontaine, J.-M. and Perrin-Riou, B. Autour des conjectures de Bloch et Kato: cohomologie galoisienne et valeurs de fonctions L. Motives (Seattle, WA, 1991), 599–706, Proc. Sympos. Pure Math., 55, Part 1, Amer. Math. Soc., Providence, RI, 1994. 

[FS] Flach, M. amd Siebel, D., Special values of the zeta-function of an algebraic surface.  2019,  arxiv: 1909.07465

[Ge] Geisser, T. Arithmetic cohomology over finite fields and special values of $\zeta$-functions, Duke Math. Journal 133 (2006), 27-57

[H] Harsthorne, R. Algebraic Geometry, Springer-Verlag (1977)

[KS] Kato, K., Saito, T., On the Conductor formula of Bloch.  Publ. Math. Inst. Hautes Etudes Sci. 100 (2004) 5-151

[Le] Levine, M. Techniques of localization in the theory of algebraic cycles, J. Alg.Geom. 10 (2001), 298-363  

[Li1] Lichtenbaum, S. The Weil-\'etale topology for number rings. Ann. of Math. (2) 170 (2009), no. 2, 657-683.

[Li2] Lichtenbaum, S.  Euler Characteristics and Special Values of Zeta-Functions, Proceedings of Fields Institute conference on algebraic cycles in honor of Spencer Bloch, Fields Institute communications, volume 56, (2009)

[Li3] Lichtenbaum, S.  Values of Zeta-Functions, Etale Cohomology, and Algebraic K-Theory. in  Springer Lecture Notes in Mathematics (342) 1973

[LR] Lichtenbaum, S. and Ramachandranm N. Values of zeta functions of algebraic surfaces at $s=1$. 2021. arxiv:2101.06530

[LS] Lichtenbaum, S, and Schlessinger, M.  The cotangent complex of a morphism. Trans. Amer. Math. Soc., 126 No. 1, (1967), 41-70

[M] Morin, B. Zeta functions of regular arithmetic schemes at s=0. Duke Math. J. 87, 1263-1336. 98/91653`---------

[R] Rapoport, M. Comparison of the regulators of Beilinson and of Borel, in Beilinson's conjectures on special values of L-functions, (M. Rapoport, N. Schappacher, P.Schneider, eds.), Perspectives in Math. 4 Acad. Press 1989

{RSS} Rapoport, M, Schappacher, N., and Schneider, S., Beilinson's conjectures on special values of L-functions, Perspectives in Math. 4 Acad, Press 1989,

[Sa] Saito, T. Parity in Bloch's conductor formula in even dimensions.  Journal de Th\'eorie des Nombres de Bordeaux 16(2004) no, 2, 403-421

[Se1] Serre, J.-P., Zeta and L-functions, in Arithmetical Algebraic Geometry, Harper and Row (1965), 82-92

[Se2] Serre, J.-P., Facteurs locaux des fonctions z\^eta des vari\'eti\'es alg\`ebriques (d\'efinitions et conjectures) S\'eminaire Delange-Pisot-Poitou, (1969/70), no. 1

[[T] Tate, J.,  On the conjectures of Birch and Swinnerton-Dyer and a geometric analog Dix expos�s sur la cohomologie des sch�mas, 189�214, Adv. Stud. Pure Math., 3, North-Holland, Amsterdam, 1968.

\bye

\bye